\documentclass[11pt]{article}  
\usepackage[greek,english]{babel}  
\usepackage[latin1]{inputenc} 
\usepackage{graphicx}
\usepackage[pdftex,bookmarks]{hyperref}
\usepackage{amsmath}
\usepackage{amssymb}
\usepackage{amsmath,amsthm}
\usepackage{verbatim}
\usepackage{sidecap}
\usepackage{tikz}

\usepackage{amsmath,amsthm}
\usepackage{cancel}
\usepackage{textcomp}
\usepackage{multicol}
\usepackage{color}
\usepackage{mathtools}  
\mathtoolsset{showonlyrefs}     

\def\ep{\epsilon}
\def\LL{{\cal L}}
\def\RR{{\mathbb R}}
\def\NN{{\mathbb N}}

\newcommand{\vs}[1]{\vskip #1pt}
\def\ni{\noindent}

\newtheorem{theorem}{Theorem}[section]

\newtheorem{lemma}[theorem]{Lemma}
\newtheorem{proposition}[theorem]{Proposition}

\theoremstyle{definition}

\begin{document}

\titlepage 
\title{\bf \Large Highly oscillatory solutions of a Neumann problem for a 
$p$-laplacian equation\thanks{Under the auspices of GNAMPA-I.N.d.A.M., Italy.
The work has been performed in the frame of the PRIN-2012-74FYK7 project
``Variational and perturbative aspects of nonlinear differential problems''.}}

\author{ALBERTO BOSCAGGIN
\quad and \quad
WALTER DAMBROSIO
}
\date{}
\maketitle

\normalsize
\begin{abstract}
We deal with a boundary value problem of the form
\begin{equation}\label{eqab}
\left\{\begin{array}{l}
\vspace{0.2cm}
-\ep(\phi_p(\ep u'))'+a(x)W'(u)=0 \\
u'(0)=0=u'(1),
\end{array}
\right.
\end{equation}
where $\phi_p(s) = \vert s \vert^{p-2} s$ for $s \in \mathbb{R}$ and $p>1$,
and $W:[-1,1] \to \RR$ is a double-well potential. We study the limit profile of solutions of \eqref{eqab} when
$\ep \to 0^+$ and, conversely, we prove the existence of nodal solutions associated with any admissible limit profile
when $\ep$ is small enough.
\end{abstract}

\noindent
{\footnotesize \textbf{AMS-Subject Classification}}. {\footnotesize 34B15, 34E15}

\section{Introduction and summary of the main results}
\def\theequation{1.\arabic{equation}}\makeatother
\setcounter{equation}{0}

As is well-known, a typical strategy to get multiplicity results for boundary value problems associated with 
nonlinear scalar second order ODEs relies on the 
investigation of the nodal properties of the solutions (see, for instance, the classical survey \cite{Ma-95}).
Quite recently, such an issue has been faced in a singular perturbation setting,
according to the following typical scheme: parameter dependent equations of the form
$$
-\ep^2 u'' + f(x,u) = 0
$$
are considered, and - for $\ep$ small enough - nodal solutions are provided, modeled 
on some limit profile for $\ep \to 0^+$ and thus exhibiting precise qualitative asymptotic properties
(depending of course on the nonlinear function $f$).
\smallbreak
In this direction, we mention on one hand the papers \cite{DeFeTa-02,FeMaTa-06,FeTo-02}, studying a one-dimensional 
Schr\"odinger equation like $-\ep^2 u'' + V(x)u - \vert u \vert^{\alpha-1} u = 0$ (with $\alpha > 1$). This line of research 
originates from the one dealing with the singularly perturbed PDE Schr\"odinger equation, which has been the object
of an enormous number of investigations in the last decades (see, among many others, \cite{AmBaCi-97,DeFe-98,Ra-92}).
On the other hand, in \cite{FeMa-03,FeMaTa-05,Na-03,NaTa-03} an equation of the type 
$-\ep^2 u'' + a(x)W'(u) = 0$, with $a$ a positive weight function and $W$ a double-well potential, is taken into account.
\smallbreak
Here, we take the work \cite{FeMaTa-05} by Felmer, Martinez and Tanaka as our starting point.
The results obtained therein, which can be applied to the spatially inhomogeneous balanced Allen-Cahn equation
\begin{equation}\label{allen}
\ep^2 u'' + a(x)u(1-u^2) = 0,
\end{equation}
and to the equation for a pendulum of variable length
\begin{equation}\label{pendu}
\ep^2 u'' + a(x)\sin(\pi u) = 0,
\end{equation}
can be roughly summarized as follows: the asymptotic behavior, for $\ep \to 0^+$, of solutions to
\eqref{allen} and \eqref{pendu} (with Neumann boundary conditions) can be characterized in term of a limit energy function
and, conversely, highly oscillatory solutions corresponding to any admissible limit profile exist for $\ep$ small
enough. More precisely, the admissible limit profiles are determined by an ordinary differential equation solved by the limit energy function and solutions to the boundary value problem are constructed using a variational approach, of broken-geodesic Nehari type 
(see also \cite{OrVe-04,TeVe-00}). Notice that this in particular shows that the above equations possess an extremely rich set of (nodal) solutions.
\smallbreak
The aim of the present paper is to extend the results in \cite{FeMaTa-05}
to equations driven by the $p$-laplacian operator. 
More precisely, throughout the paper we deal with the Neumann boundary value problem
\begin{equation}\label{bvp-1}
\left\{\begin{array}{l}
\vspace{0.2cm}
-\ep(\phi_p(\ep u'))'+a(x)W'(u)=0 \\
u'(0)=0=u'(1),
\end{array}
\right.
\end{equation}
where $\ep > 0$ and $\phi_p:\RR\to \RR$ is defined, for $p > 1$, by 
\[
\phi_p(s)=|s|^{p-2}s, \quad \forall \ s\in \RR. 
\]
As for the nonlinear term, we assume that 
$a \in C^1([0,1])$ is such that $a(x) > 0$ for every $x \in [0,1]$
and $W: [-1,1] \to \mathbb{R}$ is a $C^1$-function satisfying the following conditions:
\begin{itemize}
\item[(W1)] there exist constants $C_{-1},C_0,C_1, W_0 > 0$ such that
$$
W(u) = \frac{C_{\pm 1}}{p} \vert u - (\pm 1) \vert^p + o(\vert u - (\pm 1) \vert^p), \quad \mbox{ for } u \to \pm 1,
$$
and
$$
W(u) = W_0 - \frac{C_0}{p} \vert u \vert^p + o(\vert u \vert^p), \quad \mbox{ for } u \to 0,
$$
\item[(W2)] the function 
$$
u \in [-1,1] \setminus \{0\} \mapsto \frac{W'(u)}{\phi_p(u)}
$$
is strictly decreasing on $[-1,0)$ and strictly increasing on $(0,1]$.
\end{itemize}
Notice that from (W1) and (W2) it follows that
$$
W(\pm 1) = W'(\pm 1) = W'(0) = 0 \quad \mbox{ and } \quad W'(u)u < 0, \quad \forall \, \vert u \vert <1, u \neq 0;
$$
hence, $W$ has exactly the three critical points $\{0,\pm 1\}$:
$\pm 1$ are minima with value $0$, and $1$ is a maximum with value $W(0) = W_0 > 0$.
Typical examples of potentials $W$ satisfying the above assumptions
are for instance $W(u) = \tfrac{1}{p^2}\left( 1 - \vert u \vert^p\right)^p$, leading to the equation
\begin{equation}\label{pallen}
\ep (\phi_p(\ep u'))' + a(x)\phi_p(u) \left( 1 - \vert u \vert^p\right)^{p-1} = 0,
\end{equation}
or $W(u) = \int_u^1 \phi_p(\sin(\pi s))\,ds$, corresponding to
\begin{equation}\label{ppendu}
\ep (\phi_p(\ep u'))' + a(x)\phi_p(\sin(\pi u)) = 0.
\end{equation}
Of course, equations \eqref{pallen} and \eqref{ppendu} are natural generalizations,
to the case $p \neq 2$, of the Allen-Cahn equation \eqref{allen} and of the 
pendulum equation \eqref{pendu}, respectively.
\smallbreak
For the reader's convenience, we collect here an informal summary of the results contained in the rest of the paper.
\medbreak
\noindent
\emph{\textbf{Summary of the results.} For a family $\{ u_\ep \}$ of solutions of \eqref{bvp-1}, define the energy function
(see \eqref{def-energia} and \eqref{def-phi})
\begin{equation}\label{energia-intro}
E_\ep(x) = -\frac{p-1}{p}\frac{\ep^p}{a(x)} \vert u'_\ep(x) \vert^p + W(u_\ep(x)).
\end{equation}
Then, the following hold true.
\begin{itemize}
\item[(I)] Up to subsequences, $E_\ep$ converges for $\ep \to 0^+$ to a $C^1$ function $E$
(see Proposition \ref{energia-limitata}); moreover (see Theorem \ref{profilo-energia}) $E$ satisfies the differential equation
\begin{equation}\label{eqenergia0}
E'(x) =\dfrac{a'(x)}{a(x)} \, K(E(x)),
\end{equation}
where $K$ is a (non-Lipschitz) function - defined in \eqref{funzione-kappa} -
measuring the averaged kinetic energy of the solutions
of the autonomous equation $-(\phi_p(u'))'+W'(u)=0$.
\item[(II)] Information about the asymptotic distribution of the zeros of $u_\ep$ can
be obtained from $E$ (see Propositions \ref{stima-zeri} and \ref{accumulazione-zeri}).
\item[(III)] Equation \eqref{eqenergia0} has many solutions
(see Proposition \ref{prop-soluzioni}) and, for
any solution $E$ of it, there is a family $\{u_\ep\}$
of solutions of \eqref{bvp-1} such that its energy $E_\ep$ converges to $E$
(see Theorems \ref{main-ex-1} and \ref{main-ex-2}).
\end{itemize}
}

Let us observe that singularly perturbed equations associated with the $p$-laplacian operator
were considered for instance in \cite{Do-05,FiFu-12,Gl-10}. However, all these contributions deal with the PDE case;
we are not aware of works studying nodal solutions of ODEs driven by the $p$-laplacian
in a singular perturbation setting.
\smallbreak
To prove our results, we follow closely the approach developed in \cite{FeMaTa-05},
suitably adapted to deal with the (non-linear) differential operator $u \mapsto -\phi_p(u')'$. 
In this direction, apart from the need for quite classical phase-plane analysis tools for $p$-laplacian equations
(see, among others, \cite{DeElMa-89,MaNjZa-95,MaZa-93,Zh-01}), we point out the use of an interpolation inequality
of Landau-Kolmogorov type \cite[Chapter 1]{MiPeFi-91} which we have not found in literature and can have some independent interest
(see Lemma \ref{landaulemma}). Finally (at a more technical level), it can be worth emphasizing that many proofs in \cite{FeMaTa-05}
take advantage of a change of variable transforming the equation $-\ep^2 u'' + a(x)W'(u) = 0$ into 
$-\ep^2 u'' - \ep^2 \tfrac{a'(x)}{a(x)}u' + W'(u) = 0$, looking 
(for $\ep$ small, up to an $x$-rescaling) like a perturbation 
of an autonomous ODE. Due to the nonlinearity of the $p$-laplacian operator (when $p \neq 2$),
such a transformation is not possible for the equation in \eqref{bvp-1}; however,
all the difficulties which arise from working directly
on the original equation can be suitably overcome.
In this way, we obtain also a slightly more transparent proof of the results
in \cite{FeMaTa-05}.
\medbreak
\noindent
\textbf{Notation.}
Let us clarify that by a solution of
\eqref{bvp-1} we mean a function $u \in C^1([0,1])$,
with $u'(0) = u'(1) = 0$, such that $\phi_p(u') \in C^1([0,1])$ and the differential equation
in \eqref{bvp-1} is satisfied for every $x \in [0,1]$. By elementary regularity considerations, this is the same as
$u \in C^1([0,1])$, with $u'(0) = u'(1) = 0$,
and 
$$
\ep \int_0^1 \phi_p(\ep u'(x))\psi'(x)\,dx + \int_0^1 a(x)W'(u(x))\psi(x)\,dx = 0, \qquad \forall \, \psi \in C^{\infty}_c(]0,1[). 
$$
As usual, here $C^{\infty}_c(]0,1[)$ denotes the space of $C^\infty$ functions having support
compactly contained in the open interval $]0,1[$\,.

Throughout the paper, we use the following notation:
$p^*$ is the conjugate exponent of $p$, i.e. $1/p + 1/p^* = 1$, and
\begin{equation} \label{def-phi}
\begin{array}{l}
{\displaystyle \Phi(s)=\int_0^s \phi_p(t)\,dt=\dfrac{1}{p}|s|^p,\quad \forall s\in \RR},
\\
\displaystyle{{\Phi}_*(s)=\int_0^s \phi_p^{-1}(t)\,dt=\dfrac{1}{p^*}|s|^{p^*},\quad \forall s\in \RR,}
\\
\displaystyle{{\cal L}(s)={\Phi}_*(\phi_p(s))=\dfrac{1}{p^*}|s|^p,\quad \forall \ s\in \RR,}
\\
\displaystyle{\mathcal{L}_+^{-1}(s) =\left(\mathcal{L}|_{\mathbb{R}^+}\right)^{-1}(s) = (p^*)^{1/p} s^{1/p}, \quad \forall s\geq 0.} 
\end{array}
\end{equation}
Finally, 
$$
\pi_p = 2(p-1)^{1/p}\int_0^1\frac{ds}{(1-s^p)^{1/p}}.
$$

\section{The autonomous equation}\label{section2}
\def\theequation{2.\arabic{equation}}\makeatother
\setcounter{equation}{0}

\vs{12}
\ni
In this preliminary section we recall some facts about autonomous equations of the form
\begin{equation} \label{autonoma}
-(\phi_p(v'))'+a W'(v)=0,
\end{equation}
where $a > 0$ is a constant. It is well known that the energy function
\[
H_a(p,q)=-\dfrac{1}{a} \LL(p)+W(q),\qquad (p,q)\in \RR^2,
\]
is preserved along the solutions of \eqref{autonoma}:
that is, if $v$ solves \eqref{autonoma} then 
$H_a(v,v') \equiv \textnormal{const}$. Moreover:
\begin{itemize}
\item $H_a(v,v') \equiv \xi \in (0,W_0)$ implies that $v$ is periodic,
\item $H_a(v,v') \equiv W_0$ implies $v \equiv 0$,
\item $H_a(v,v') \equiv 0$ implies either $v \equiv \pm1$, or $v$ is an heteroclinic
solution joining the rest points $-1$ and $1$ (see the proof of Proposition \ref{eterocline}).
\end{itemize}

We first prove that solutions with zero energy are of the form specified above and that their kinetic energy is 
in $L^1(\RR)$.
\vs{12}
\ni
\begin{proposition} \label{eterocline} 
Let us consider a solution $v$ of \eqref{autonoma} with energy $H_a(v,v') \equiv 0$. Then $v$ is defined in $\RR$ and
\begin{equation} \label{energia-integrabile}
\int_{-\infty}^{+\infty} \LL (v'(x))\,dx<+\infty.
\end{equation}
\end{proposition}
\vs{8}
\ni
\begin{proof}
Clearly the result is true for the equilibrium solutions $v \equiv -1$ and $v \equiv 1$.
Hence, we can assume that $v$ is non-constant; without loss of generality, we also suppose
that $v(0) = 0$ and $v'(0) > 0$, and we consider $x \geq 0$ 
(the arguments for the other cases are analogous). We observe 
that $v$ is continuable in the future as long as $v\neq 1$; assume then that
\[
\lim_{x\to \Lambda^-} v(x)=1,
\]
where $\Lambda>0$ denotes the right extremum of the maximal interval of existence of $v$. We prove that indeed $\Lambda = +\infty$.

Let us write the conservation of energy
\begin{equation} \label{d11}
-\dfrac{1}{a} \mathcal{L}(v'(x))+W(v(x))=0,\quad \forall \ x\in (0,\Lambda);
\end{equation}
from assumption (W1) we know that
\begin{equation} \label{d12}
W(u)=\widehat C_1(1-u)^p+R_1(u),
\end{equation}
with $\widehat C_1 = C_1/p$ and $R_1(u)=o((1-u)^p)$, $u\to 1$. From \eqref{d11} and \eqref{d12} we deduce that
\begin{equation} \label{d13}
v'(x)=(1-v(x))\sqrt[p]{a p^* \widehat C_1+\dfrac{R_1(v(x))}{(1-v(x))^p}},\quad \forall \ x\in (0,\Lambda);
\end{equation}
as a consequence, we infer that there exists $\Lambda_\infty>0$ such that
\begin{equation} \label{d14}
\dfrac{1}{2}\ \sqrt[p]{a p^* \widehat C_1}\leq \dfrac{v'(x)}{1-v(x)}\leq \dfrac{3}{2}\ \sqrt[p]{ap^* \widehat C_1},\quad \forall \ x\in (\Lambda_\infty ,\Lambda).
\end{equation}
By integrating, we obtain
\begin{equation} \label{d15}
\dfrac{L}{2}(x-\Lambda_\infty)-K\leq -\log (1-v(x))\leq \dfrac{3}{2}L(x-\Lambda_\infty)-K,\quad \forall \ x\in (\Lambda_\infty ,\Lambda),
\end{equation}
for some constants $L>0$ and $K\in \RR$. Passing to the limit in \eqref{d15} for $x\to \Lambda^-$, we plainly deduce that $\Lambda =+\infty$; taking the exponentials, we also obtain
\begin{equation} \label{d16}
K'\ e^{-\dfrac{3}{2}L(x-\Lambda_\infty)}\leq 1-v(x)\leq K'\ e^{-\dfrac{L}{2}(x-\Lambda_\infty)},\quad \forall \ x>\Lambda_\infty,
\end{equation}
for some $K'>0$. This proves that $v$ goes exponentially to $1$ as $x\to +\infty$; this is sufficient to show also that
\[
\int_0^{+\infty} \LL (v'(x))\,dx<+\infty.
\]
Indeed, from the conservation of the energy we have
\[
\LL (v'(x))=a W(v(x))\sim a \widehat C_1 (1-v(x))^p,\quad x\to +\infty,
\]
which is integrable since \eqref{d16} holds true.
\end{proof}
\vs{12}
\ni
In what follows, let $\bar{v}_{a,\xi}(\cdot)$ denote:
\begin{itemize}
\item if $\xi \in (0,W_0)$, the (unique) periodic solution of \eqref{autonoma}
with energy $\xi$, and such that $\bar{v}_{a,\xi}(0) = 0$,
$\bar{v}'_{a,\xi}(0) > 0$ (of course, any other solutions having energy $\xi$ is a translation
of $\bar{v}_{a,\xi}$);
\item if $\xi = W_0$, the trivial solution $\bar{v}_{a,\xi}\equiv 0$;
\item if $\xi = 0$, the heteroclinic solution with $\bar{v}_{a,\xi}(0) = 0$,
$\bar{v}'_{a,\xi}(0) > 0$
(any other non-constant solution with zero energy is a translation
of $\bar{v}_{a,\xi}(x)$ or of $\bar{v}_{a,\xi}(-x)$).
\end{itemize}
For $\xi\in (0,W_0)$, we also denote by $T_{a}(\xi)$ be the period of $\bar{v}_{a,\xi}$;
as it is well-known, this can be computed via the time-map formula
\begin{equation} \label{periodo}
T_{a}(\xi)=2\int_{h_-(\xi)}^{h_+(\xi)} \dfrac{ds}{\LL_+^{-1} (a(W(s)-\xi))},
\end{equation}
where $h_-(\xi)<0<h_+(\xi)$ are the unique points such that $W(h_\pm (\xi))=\xi$.
Notice that here we have already used the assumptions (W1) and (W2).
\vs{6}
\ni
Finally, we define the averaged kinetic energy of $\bar{v}_{a,\xi}$, namely
\begin{equation} \label{funzione-kappa}
K_a(\xi)=\dfrac{1}{T_a(\xi)} \int_0^{T_a(\xi)} \LL (\bar{v}_{a,\xi}'(y))\,dy.
\end{equation}
In the particular case $a=1$ we use the notation $T$ and $K$.
\vs{12}
\ni
\begin{lemma} \label{confr-periodi-kappa}
For every $\xi\in (0,W_0)$ we have
\begin{equation} \label{relaz-periodi-kappa}
T(\xi)=a^{1/p} T_a(\xi) \qquad \mbox{ and } \qquad
K(\xi)=a^{-1} K_a(\xi).
\end{equation}
\end{lemma}
\vs{8}
\ni
\begin{proof}
Let us first observe that if $v$ is a solution of \eqref{autonoma} then the function $v_*:\RR\to \RR$ defined by
\begin{equation} \label{aa1}
v_*(y)=v\left(\dfrac{y}{a
^{1/p}}\right),\quad \forall \ y\in \RR,
\end{equation}
is a solution of \eqref{autonoma} with $a=1$; we denote by $\xi$ and $\xi'$ the levels of energy of these solutions, respectively. From \eqref{aa1} we infer that
\[
T_a(\xi)=\dfrac{T(\xi')}{a^{1/p}}.
\]
Since
\begin{align*}
H_a(v,v') & =-\dfrac{1}{a} \LL (v')+W(v)=-\dfrac{1}{a} \LL (a^{1/p} v'_*)+W(v_*) \\
& =-\LL (v'_*) *W(v_*)=H_1(v_*,v'_*),
\end{align*}
we deduce that $\xi=\xi'$; this completes the proof of the first relation in \eqref{relaz-periodi-kappa}. 

Finally, a simple computation shows that
\begin{align*}
K_a(\xi) & =\dfrac{1}{T_a(\xi)}\ \int_0^{T_a(\xi)} \LL(v'_0(y))\,dy= \dfrac{1}{T_a(\xi)}\ \int_0^{T_a(\xi)} \LL(a^{1/p} v'_* (a^{1/p}y))\,dy \\
& =\dfrac{a}{T_a(\xi)}\ \int_0^{T_a(\xi)} \LL(v'_* (a^{1/p}y))\,dy=\dfrac{a}{a^{1/p}T_a(\xi)}\ \int_0^{a^{1/p}T_a(\xi)} \LL(v'_* (u))\,du \\
&=\dfrac{a}{T(\xi)}\ \int_0^{T(\xi)} \LL(v'_* (u))\,du=a K(\xi),
\end{align*}
proving the second relation in \eqref{relaz-periodi-kappa}.
\end{proof}
\vs{12}
\ni
The next propositions collect some properties of the functions $T$ e $K$
which are needed in the rest of the paper.
\vs{12}
\ni
\begin{proposition} \label{prop-periodo} The function $T$ satisfies:
\begin{itemize}
\item[(i)] $T\in C^1((0,W_0))$ and $T'<0$ in $(0,W_0)$;
\item[(ii)] $T(\xi)\to \dfrac{2\pi_p}{\sqrt[p]{C_0}}$, for $\xi \to W_0^-$;
\item[(iii)] $T(\xi) \to +\infty$ for $\xi \to 0^+$ and, more precisely, 
$T(\xi) \sim -C^*\log \xi$, for some $C^*>0$.
\end{itemize}
\end{proposition}
\vs{8}
\ni
\begin{proof}
We first observe that (i) can be proved as in the case $p=2$ (see for instance \cite{Op-61}), while the proof of (ii) can be found in \cite{MaZa-93}; we give the proof of (iii) only. We write the conservation of energy as
\[
\LL (v')+F(v)=(W_0-\xi),
\]
where $F(s)=(W_0-W(s))$, for every $s\in (-1,1)$,
and $v = \bar{v}_{1,\xi}$. Without loss of generality we can assume that 
\[
(v(0),v'(0))=(0,\LL^{-1}_+(W_0-\xi)),\quad (v(T/4),v'(T/4))=(F^{-1}_+ (W_0-\xi),0);
\] 
let us denote $\alpha (\xi)=F^{-1}_+ (W_0-\xi)$.
We then have
\begin{equation} \label{c11}
\begin{array}{l}
\displaystyle{4T(\xi)=\int_0^{\alpha(\xi)} \dfrac{du}{\LL^{-1}_+(F(\alpha (\xi)-F(u))}}\\
\\
=\displaystyle{\alpha (\xi) \int_0^1 \dfrac{dt}{\LL^{-1}_+(F(\alpha (\xi)-F(\alpha (\xi)t))}=\dfrac{\alpha (\xi)}{(p^*)^{1/p}} \int_0^1 \dfrac{dt}{\sqrt[p]{W(\alpha t)-W(\alpha)}}}\\
\\
=\displaystyle{\dfrac{\alpha (\xi)}{(p^*)^{1/p}} \left(\int_0^\delta \dfrac{dt}{\sqrt[p]{W(\alpha t)-W(\alpha)}}+\int_\delta^1 \dfrac{dt}{\sqrt[p]{W(\alpha t)-W(\alpha)}}\right),}
\end{array}
\end{equation}
for every $\delta\in (0,1)$ to be chosen.
\ni
Now, let us observe that $\xi \to 0^+$ implies $\alpha (\xi)\to 1$; more precisely, using 
(W1) it is possible to show that there exists $C'>0$ such that
\begin{equation} \label{cc1}
\alpha (\xi)=1-C' \xi^{1/p} +o(\xi^{1/p}),\qquad \xi \to 0^+.
\end{equation}
Noting also that
\[
\lim_{\alpha \to 1} \int_0^\delta \dfrac{dt}{\sqrt[p]{W(\alpha t)-W(\alpha)}}<+\infty,
\]
for every $\delta >0$, equation \eqref{c11} can be written as
\begin{equation} \label{cc2}
4T(\xi)=O(1)+\dfrac{\alpha (\xi)}{(p^*)^{1/p}} \int_\delta^1 \dfrac{dt}{\sqrt[p]{W(\alpha t)-W(\alpha)}},\qquad \xi \to 0^+.
\end{equation}
We now show that
\begin{equation} \label{c12}
\int_\delta^1 \dfrac{dt}{\sqrt[p]{W(\alpha t)-W(\alpha)}} \sim -\log (1-\alpha),\qquad \alpha \to 1.
\end{equation}
From assumption (W1) we know that 
\begin{equation} \label{c13}
W(\alpha t)-W(\alpha)=C_1 (|\alpha t-1|^p-|\alpha -1|^p)+R_1(\alpha)+R_2(\alpha t),
\end{equation}
with $R_1(t)=o(|\alpha -1|)$, $\alpha \to 1$ and $R_2(\alpha t)=o(|\alpha t-1|)$, $\alpha t\to 1$. Hence, when $\alpha$ and $\delta$ are close to $1$, there exist two constants $A$ and $B$ such that
\[
B^{-p} (|\alpha t-1|^p-|\alpha -1|^p)\leq W(\alpha t)-W(\alpha)\leq A^{-p} (|\alpha t-1|^p-|\alpha -1|^p);
\]
as a consequence, for suitable $\alpha$ and $\delta$, we have
\begin{equation} \label{c14}
A \int_\delta^1 \dfrac{dt}{\sqrt[p]{|\alpha t-1|^p-|\alpha -1|^p}}\leq \int_\delta^1 \dfrac{dt}{\sqrt[p]{W(\alpha t)-W(\alpha)}} \leq B \int_\delta^1 \dfrac{dt}{\sqrt[p]{|\alpha t-1|^p-|\alpha -1|^p}}.
\end{equation}
Now, it is possible to show that there are constants $M_p$ and $N_p$ such that
\begin{equation} \label{c15}
0<M_p |x-y|\ (|x|+|y|)^{p-1}\leq ||x|^{p-1}x-|y|^{p-1}y|\leq N_p |x-y|\ (|x|+|y|)^{p-1},\quad \forall \ (x,y)\in \RR^2.
\end{equation}
Therefore, from \eqref{c14} and \eqref{c15} we infer that it is sufficient to study 
\begin{equation} \label{c16}
\int_\delta^1 \dfrac{dt}{\alpha \sqrt[p]{(1-t)(2/\alpha -(t+1))^{p-1}}}
\end{equation}
when $\alpha \to 1$. A first change of variables leads to
\begin{equation} \label{c17}
\int_0^{1-\delta} \dfrac{ds}{\sqrt[p]{s(s+2/\alpha-2)^{p-1}}}=\int_0^{1-\delta} \dfrac{ds}{\sqrt[p]{s(s+\lambda)^{p-1}}},
\end{equation}
with $\lambda=2/\alpha-2$; with the change of variables $s=\lambda r$ the integral in \eqref{c17} reduces to
\begin{equation} \label{c18}
\int_0^{(1-\delta)/\lambda} \dfrac{dr}{\sqrt[p]{r(r+1)^{p-1}}}.
\end{equation}
Finally, we note that $\lambda \to 0^+$ when $\alpha \to 1^-$; moreover, we have
\begin{equation} \label{c19}
\int_0^{(1-\delta)/\lambda} \dfrac{dr}{\sqrt[p]{r(r+1)^{p-1}}} \sim \log \dfrac{1-\delta}{\lambda},\qquad \lambda \to 0^+.
\end{equation}
From \eqref{c19} we obtain that
\begin{equation} \label{c20}
\int_\delta^1 \dfrac{dt}{\alpha \sqrt[p]{(1-t)(2/\alpha -(t+1))^{p-1}}} \sim \log \dfrac{1-\delta}{2}-\log (1-\alpha),\qquad \alpha \to 1;
\end{equation}
from \eqref{c14} and \eqref{c20} we deduce the validity of \eqref{c12}.

Now, from \eqref{cc1} we infer that
\begin{equation} \label{cc3}
\log (1-\alpha)\sim \dfrac{1}{p} \log \xi,\quad \xi \to 0^+.
\end{equation}
From \eqref{cc2}-\eqref{c12} and \eqref{cc3} we obtain the thesis.
\end{proof}
\vs{12}
\ni
\begin{proposition} \label{prop-kappa} 
The function $K$ satisfies:
\begin{itemize}
\item[(i)] $K\in C^1((0,W_0))$ and $K > 0$ in $(0,W_0)$
\item[(ii)]$K(\xi) \to 0$ both for $\xi \to 0^+$ and for $\xi \to W_0^-$;
\item[(iii)] $\displaystyle{\int_0^{W_0/2} \dfrac{d\xi}{K(\xi)}} < +\infty$;
\item[(iv)] $\displaystyle{\int_{W_0/2}^{W_0} \dfrac{d\xi}{K(\xi)}=+\infty}$.
\end{itemize}
\end{proposition}
\vs{8}
\ni
\begin{proof}
The proof can be obtained by arguing as in \cite[Lemma 3.2]{FeMaTa-05}, taking into account
Proposition \ref{eterocline}.
\end{proof}
\vs{12}
\ni
We end this section with two more technical results; the second one will be used in the proof of Lemma \ref{stimeminimi}.
\vs{12}
\ni
\begin{proposition} \label{Nakashima-1} There exists $M_0>0$ such that for every $M>M_0$ there exist a unique positive solution $v_+$ and a unique negative solution $v_-$ of \eqref{autonoma} such that $v_\pm (-M)=0=v_\pm (M)$. Moreover, there exist $R_\pm>0$ and $C_\pm >0$ such that for every $M>M_0$ we have
\begin{equation} \label{stime-asintotiche-Nakashima}
\begin{array}{l}
C_-e^{-R_- M}\leq 1-v_+(0)\leq C_+e^{-R_+ M}\\
\\
C_-e^{-R_- M}\leq v_-(0)+1\leq C_+e^{-R_+ M}.
\end{array}
\end{equation}
\end{proposition}
\vs{8}
\ni
\begin{proof}
The existence and uniqueness of $v_\pm$ plainly follow from the monotonicity properties of the time map $T$; using the notation of the proof of Proposition \ref{prop-periodo}, we have $T=2M$ and $v_\pm (0)=\pm \alpha=\pm F^{-1}_+(W_0-\xi)$. 
\vs{4}
\ni
We prove the first inequality in \eqref{stime-asintotiche-Nakashima}; the proof of the second one is analogous. First of all, let us observe that when $M\to +\infty$, then $\xi \to 0^+$; as a consequence, $\alpha=v_+(0)\to 1^-$.
\ni
Moreover, from \eqref{cc2}-\eqref{c12} we know that there exist $K_\pm>0$ and $C>0$ such that
\[
K_-\leq \dfrac{2M-2C}{-\log (1-v_+(0))}\leq K_+,\quad v_+(0)\to 1^-.
\] 
The result follows by solving these inequalities with respect to $1-v_+(0)$.
\end{proof}
\vs{12}
\ni
Arguing as in \cite[Prop. 2.4-2.5]{Na-03}, from Proposition \ref{Nakashima-1} we deduce the following result (see also \cite[Lemma A.11]{FeMaTa-05}):
\vs{12}
\ni
\begin{proposition} \label{Nakashima-2} There exist $K_1, K_2>0$ such that,
for every $\ep > 0$, for every $[s,t]\subset [0,1]$ and for every solution $u$ of the equation in \eqref{bvp-1} 
on $[s,t]$ and with $u(x)\in [0,1]$ for $x \in [s,t]$, 
\begin{equation} \label{stima-asintotica-Nakashima-2}
\begin{array}{l}
|u(x)-1|+|\ep u'(x)|\leq K_1 e^{-K_2 \min (|x-s|,|x-t|)/\ep},\quad \forall \ x\in [s,t].
\end{array}
\end{equation}
\end{proposition}
\vs{12}
\ni
\section{The limit energy function}\label{section3}
\def\theequation{3.\arabic{equation}}\makeatother
\setcounter{equation}{0}
\vs{12}
\ni
In this section, we define an energy function for solutions of the non-autonomous problem
\eqref{bvp-1} (compare with \eqref{energia-intro}) and we study its convergence to a limit profile.
\medbreak
\ni
We start with a simple lemma.
For its proof, we need to observe that, for any $v \in C^1(\RR)$
such that $\phi_p(v') \in C^1(\RR)$, the following elementary inequality holds true:
\begin{equation}\label{interpolazionefacile}
\Vert \phi_p(v') \Vert_\infty \leq 2^{p-1} \Vert v \Vert_\infty^{p-1}  + \Vert (\phi_p(v'))' \Vert_\infty.
\end{equation}
\vs{12}
\ni
\begin{lemma}\label{conv-riscalate0}
Let $\{u_\ep\}$ be a family of solutions of \eqref{bvp-1}, fix $x_0\in (0,1)$ and define
\begin{equation} \label{def-riscalate0}
v_{x_0,\ep}(y)=u_\ep(x_0+\ep y),\quad \forall \ y \in \left[\dfrac{-x_0}{\ep},\dfrac{1-x_0}{\ep} \right]. 
\end{equation}
Then, up to subsequences, for $\ep \to 0^+$ the family $\{v_\ep\}$ converges in $C^1_{\textnormal{loc}}$ to a function
$\bar{v} \in C^1(\RR)$ solving the differential equation
\begin{equation}\label{eq-limite0}
-(\phi_p(\bar{v}'))'+a(x_0)W(\bar{v})=0\quad \mbox{ in } \; \RR.
\end{equation}
\end{lemma}
\vs{8}
\ni
\begin{proof}
We first claim that $v_{x_0,\ep}$ is a solution of the differential equation
\begin{equation} \label{2-222}
-(\phi_p(v_{x_0,\ep}'(y)))'+a(x_0+\ep y)W(v_{x_0,\ep}(y))=0\quad \mbox{in}\ \left[ -\dfrac{x_0}{\ep},\dfrac{1-x_0}{\ep}\right].
\end{equation}
This is almost obvious, but we give the details.
By definition, \eqref{2-222} means
\begin{equation} \label{2-333}
\displaystyle{\int_{-x_0/\ep}^{(1-x_0)/\ep} \phi_p(v_{x_0,\ep}'(y))\psi'(y)\,dy=\int_{-x_0/\ep}^{(1-x_0)/\ep} a(x_0+\ep y)
W(v_{x_0,\ep}(y))\psi(y)\,dy,}
\end{equation}
for every $\psi\in C^\infty_0\left( \left] -{x_0}/{\ep},{1-x_0}/{\ep}\right[ \right)$.
Define $\widetilde\psi \in C^\infty_0(]0,1[)$ by setting
$$
\widetilde\psi(x) = \psi\left( \frac{x-x_0}{\ep}\right), \quad \forall \ x\in (0,1);
$$
since $u_\ep$ satisfies the differential equation in \eqref{bvp-1},
we have
$$
\int_{0}^{1} \phi_p(\ep u'_\ep(x))\widetilde\psi'(x)\,dx=\frac{1}{\ep}\int_{0}^{1} a(x)W(u_\ep(x))\widetilde\psi(x)\,dx.
$$
Changing variable in the integral, we thus get \eqref{2-222}.

\ni
Since $\Vert v_{x_0,\ep} \Vert_\infty \leq 1$, the differential equation
\eqref{2-222} implies (via the inequality \eqref{interpolazionefacile}) that
$v_{x_0,\ep}$ is bounded in the $C^1$-norm. Hence Ascoli-Arzela theorem
ensures that $v_{x_0,\ep} \to \bar{v}$ in $C_{\textnormal{loc}}$, for a suitable $\bar{v} \in C(\mathbb{R})$.
Now, a standard argument (using again \eqref{interpolazionefacile})
permits to pass to the limit in \eqref{2-222} so that $\bar{v} \in C^1(\RR)$ and solves  
\eqref{eq-limite0}.
\end{proof}
\vs{12}
\ni
To proceed further, we need to establish a Landau-Kolmogorov inequality
for the $\phi_p$-operator.
\vs{12}
\ni
\begin{lemma}\label{landaulemma}
For every $v \in C^1([0,1])$ such that
$v'(0) = v'(1) = 0$ and $\phi_p(v') \in C^1([0,1])$, it holds
\begin{equation}\label{landau}
\Vert \phi_p(v') \Vert_{\infty}^{\frac{p}{p-1}} \leq 4 \gamma_p \Vert v \Vert_{\infty} \Vert (\phi_p(v'))' \Vert_{\infty},
\end{equation}
where $\gamma_p = 1/(p-1)$ if $1 < p < 2$ and $\gamma_p = 1$ if $p \geq 2$.
\end{lemma}
\vs{8}
\ni
\begin{proof}
Let us select $x_1 \in [0,1]$ such that $\vert \phi_p(v'(x_1))) \vert = \Vert \phi_p(v') \Vert_\infty$
and assume, w.l.o.g, that $\phi_p(v'(x_1)) > 0$. Also, let $x_2 \in \,]x_1,1]$ be
the first point such that $\phi_p(v'(x_2)) = 0$ (which of course exists, since 
$\phi_p(v'(1)) = 0$) and define the function
$w: [x_1,+\infty)$ by setting 
$$
\begin{array}{ll}
\vspace{0.2cm}
\displaystyle{w(x) = \int_{x_1}^x \phi_p(v'(s))^{\gamma_p}\,ds 
= \int_{x_1}^x v'(s)^{\gamma_p (p-1)}}, & \quad 
\mbox{ for } x \in [x_1,x_2] \\
w(x) = w(2x_2 - x), &\quad \mbox{ for } x \in [x_2,2x_2 - x_1] 
\end{array}
$$
and then extending by $2(x_2-x_1)$-periodicity. Notice that
$w \in W^{2,\infty}(x_1,+\infty)$, with
\begin{equation}\label{derivatew}
\vert w'(x) \vert = \vert \phi_p(v') \vert^{\gamma_p} 
\quad \mbox{ and } \quad
\vert w''(x) \vert = \gamma_p \vert \phi_p(v') \vert^{\gamma_p-1} \vert (\phi_p(v'))' \vert; 
\end{equation}
moreover
\begin{equation}\label{stimaw}
\Vert w \Vert_\infty = w(x_2) \leq \Vert v' \Vert_{\infty}^{\gamma_p(p-1)-1}\int_{x_1}^{x_2} v'(s)\,ds \leq 2 \Vert \phi_p(v') \Vert_{\infty}^
{\gamma_p-\tfrac{1}{p-1}} \Vert v \Vert_\infty.
\end{equation}
Consider now the function 
$$
\widetilde{w}(x) = w(x) + \frac{1}{2}\gamma_p
\Vert \phi_p(v') \Vert_\infty^{\gamma_p-1} \Vert (\phi_p(v'))' \Vert_\infty
(x-x_1)^2, \quad x \in [x_1,+\infty);
$$
using \eqref{derivatew}, we easily see that $\widetilde{w}'' \geq 0$ a.e..
Hence, for $x \geq x_1$,
$$
\widetilde w'(x_1)(x-x_1) \leq \widetilde w(x) - \widetilde w(x_1) = w(x),
$$
which implies, taking into account the definition of $w$
and \eqref{stimaw}
$$
\Vert \phi_p(v') \Vert_\infty^{\gamma_p} \leq \frac{ 2 \Vert \phi_p(v') \Vert_{\infty}^
{\gamma_p-\tfrac{1}{p-1}} \Vert v \Vert_\infty}{x - x_1} + \frac{1}{2}\gamma_p
\Vert \phi_p(v') \Vert_\infty^{\gamma_p-1} \Vert (\phi_p(v'))' \Vert_\infty(x-x_1). 
$$
Minimizing the above expression on $[x_1,+\infty)$, we obtain
\begin{align*}
\Vert \phi_p(v') \Vert_\infty^{\gamma_p} & \leq 2 \gamma_p^{1/2}
\Vert \phi_p(v') \Vert_{\infty}^{\left(\tfrac{\gamma_p}{2} - \tfrac{1}{2(p-1)} + \tfrac{\gamma_p-1}{2}\right)} \Vert v \Vert_\infty^{1/2}
\Vert (\phi_p(v'))' \Vert_{\infty}^{1/2} \\
& = 2 \gamma_p^{1/2}\Vert \phi_p(v') \Vert_{\infty}^{\gamma_p - \tfrac{p}{2(p-1)}} \Vert v \Vert_\infty^{1/2}
\Vert (\phi_p(v'))' \Vert_{\infty}^{1/2},
\end{align*}
thus concluding the proof.
\end{proof}
\vs{12}
\ni
\subsection{A priori estimates and convergence results}
\vs{12}
\ni
Let $\{u_\ep\}$ be a family of solutions of \eqref{bvp-1} and define the energy
\begin{equation} \label{def-energia}
E_\ep (x)=-\dfrac{1}{a(x)} \mathcal{L}(\ep u'_\ep(x)) +W(u_\ep (x)),\quad \forall \ x\in [0,1].
\end{equation}
\vs{12}
\ni
\begin{proposition}\label{energia-limitata}
The energy satisfies the differential equation
\begin{equation}\label{eq-energia}
E'_\ep (x)=\dfrac{a'(x)}{a(x)^2}\,\mathcal{L}(\ep u'_\ep(x)),\quad \forall \ x\in [0,1].
\end{equation}
Moroever, up to subsequence,
\begin{equation} \label{conv-energia}
E_\ep\to E,\quad {\mbox{uniformly in $[0,1]$.}}
\end{equation}
\end{proposition}
\vs{8}
\ni
\begin{proof}
Recalling the definition of $\mathcal{L}$,
for every $x\in [0,1]$ we have
\begin{equation} \label{111}
E'_\ep (x) = \dfrac{a'(x)}{a(x)^2} \mathcal{L}(\ep u'_\ep(x)) -\dfrac{1}{a(x)} {\Phi'}_* (\phi_p(\ep u'_\ep(x)) 
(\phi_p(\ep u'_\ep (x)))'+W'(u_\ep (x))u'_\ep (x).
\end{equation}
Taking into account that ${\Phi'}_* (\phi_p(s)) = s$ and the differential equation,
$$
E'_\ep (x) = \dfrac{a'(x)}{a(x)^2} \mathcal{L}(\ep u'_\ep(x))- \dfrac{1}{a(x)} \left(\ep u'_\ep (x)\dfrac{1}{\ep} a(x)W'(u_\ep(x))-a(x) W'(u_\ep (x))u'_\ep (x)\right),
$$
whence \eqref{eq-energia} follows.
\vs{4}
\ni
As for the convergence, we first observe that 
$$
||W(u_\ep)||_\infty \leq W_0,\quad \forall \ep > 0,
$$
so that 
$$
||\ep (\phi_p(\ep u'_\ep))'||_\infty \leq \Vert a \Vert_\infty W_0,\quad \forall \ep > 0.
$$
Using the Landau-Kolmogorov inequality \eqref{landau} with $v = \ep u'_\ep$, we deduce that there exists $M>0$ such that
\begin{equation} \label{uprimo-lim}
||\phi_p (\ep u'_\ep)||_\infty\leq M,\quad \forall \ep > 0.
\end{equation}
Recalling \eqref{def-energia} and \eqref{eq-energia}, this proves the uniform boundedness of $E_\ep$ in $W^{1,\infty} (0,1)$,
from which \eqref{conv-energia} follows.
\end{proof}
\vs{12}
\ni
Notice that, combining Lemma \ref{conv-riscalate0} with Proposition \ref{energia-limitata},
we can now state the following result.
\vs{12}
\ni
\begin{proposition} \label{conv-riscalate}
Let $\{v_{x_0,\ep}\}$ be 
as in \eqref{def-riscalate0}.
Then, up to subsequences, $v_{x_0,\ep}$ converges in $C^1_{\textnormal{loc}}$ to a function $\bar{v}$ such that
\begin{equation} \label{eq-limite}
\left\lbrace \begin{array}{l}
-(\phi_p(\bar{v}'))'+a(x_0)W(\bar{v})=0\quad \mbox{in} \ \RR\\
\\
-\dfrac{1}{
a(x_0)}\mathcal{L}(\bar{v}')+W(\bar{v})\equiv E(x_0),
\end{array}
\right.
\end{equation}
where $E$ is the limit profile of the energy, given in \eqref{conv-energia}.
Hence, $E(x_0) \in [0,W_0]$ and, according to the notation of Section \ref{section2}:
\begin{itemize}
\item if $E(x_0) \in (0,W_0)$, then $v_{x_0,\ep}(x) \to \bar{v}_{a(x_0),E(x_0)}(x + t(x_0))$
for some $t(x_0) \in \mathbb{R}$;
\item if $E(x_0) = W_0$, then $v_{x_0,\ep} \to 0$;
\item if $E(x_0) = 0$, then either $v_{x_0,\ep} \to \pm 1$
or $v_{x_0,\ep}(x) \to \bar{v}_{a(x_0),E(x_0)}(\pm x + t(x_0))$ for some $t(x_0) \in \mathbb{R}$.
\end{itemize}
\end{proposition}
\vs{12}
\ni
\subsection{The limit equation}
\vs{12}
\ni
In this section we prove that the function $E$ given in \eqref{conv-energia}
satisfies a first order differential equation.
\vs{12}
\ni
\begin{theorem} \label{profilo-energia}
The function $E$ satisfies the differential equation
\begin{equation} \label{eq-profilo}
E'(x)=\dfrac{a'(x)}{a(x)} \, K(E(x)),
\end{equation}
where the function $K$ is defined in \eqref{funzione-kappa}.
\end{theorem}
\vs{8}
\ni
\begin{proof} Let us write 
\begin{equation} \label{energia-debole}
\int_0^1 E'_\ep (x)\psi(x)\,dx=-\int_0^1 E_\ep (x)\psi'(x)\,dx,
\end{equation}
for every $\ep >0$ and for every $\psi \in C^{\infty}_0(]0,1[).$ Up to a subsequence, from \eqref{conv-energia} we plainly deduce that
\[
\lim_{\ep \to 0^+} \int_0^1 E_\ep (x)\psi'(x)\,dx=\int_0^1 E(x)\psi'(x)\,dx;
\]
on the other hand, from \eqref{eq-energia} we have
\[
\lim_{\ep \to 0^+} \int_0^1 E'_\ep (x)\psi(x)\,dx=\lim_{\ep \to 0^+} \int_0^1 \dfrac{a'(x)}{a(x)^2}\mathcal{L}(\ep u'_\ep (x))\psi(x)\,dx.
\]
Hence, we have to prove that
\begin{equation} \label{n1}
\lim_{\ep \to 0^+} \int_0^1 \dfrac{a'(x)}{a(x)^2}\mathcal{L}(\ep u'_\ep (x))\psi(x)\,dx=\int_0^1 \dfrac{a'(x)}{a(x)}K(E(x))\psi(x)\,dx, 
\end{equation}
for every $\psi \in C^{\infty}_0(]0,1[).$
\vs{4}
\ni
We argue as follows. For every $s>0$ let $\rho_s:\RR\to \RR$ be defined as
\[
\rho_s(x)=
\left\{
\begin{array}{ll}
1/s&{\mbox{if $x\in [-s/2,s/2]$}}\\
&\\
0&{\mbox{otherwise}}
\end{array}
\right.
\]
and let, for every $L>0$, $\psi_{\ep,L}=\rho_{\ep L} \ast \psi$, where $\ast$ denotes the convolution
product. Since, by standard properties of convolution, for any fixed $L > 0$ it holds,
\[
\lim_{\ep \to 0^+} ||\psi-\psi_{\ep,L}||_{L^{\infty}(0,1)}=0,
\]
an elementary dominated convergence argument (based on \eqref{uprimo-lim} too) shows that,
for every $L > 0$,
\[
\lim_{\ep \to 0^+} \int_0^1 \dfrac{a'(x)}{a(x)^2}\mathcal{L}(\ep u'_\ep (x))\left(\psi(x)-\psi_{\ep L}(x)\right)\,dx=0.
\]
Hence, we can prove \eqref{n1} by showing that
\begin{equation} \label{n2}
\lim_{\ep \to 0^+} \int_0^1 \dfrac{a'(x)}{a(x)^2}\mathcal{L}(\ep u'_\ep (x))\psi_{\ep L}(x)\,dx=\int_0^1 \dfrac{a'(x)}{a(x)}K(E(x))\psi(x)\,dx, 
\end{equation}
for every $L>0$.
\vs{4}
\ni
Using standard properties of the convolution, we can write
\begin{equation} \label{n3}
\int_0^1 \dfrac{a'(x)}{a(x)^2}\mathcal{L}(\ep u'_\ep (x)))\psi_{\ep L}(x)\,dx=
\int_0^1 F_{\ep,L}(x)\psi(x)\,dx,
\end{equation}
where
\begin{align*}
F_{\ep,L}(x) & =
\left(\rho_{\ep L} \ast \dfrac{a'(\cdot )}{a(\cdot )^2}\mathcal{L}(\ep u'_\ep (\cdot))\right)(x) \\
& =\int_{\RR} \rho_{\ep L} (y) \dfrac{a'(x-y )}{a(x-y )^2}\mathcal{L}(\ep u'_\ep (x-y))\,dy \\
&=\dfrac{1}{L} \int_{-L/2}^{L/2} \dfrac{a'(x+\ep z )}{a(x+\ep z )^2}\mathcal{L}(\ep u'_\ep (x+\ep z))\,dz\\
&=\dfrac{1}{L} \int_{-L/2}^{L/2} \dfrac{a'(x+\ep z )}{a(x+\ep z )^2}\mathcal{L}(v'_{x,\ep} (z))\,dz,
\end{align*}
and $v_{x,\ep}$ is defined as in \eqref{def-riscalate0},
that is, $v_{x,\ep}(z) = u_\ep(x + \ep z)$. 
We thus observe that, from Proposition \ref{conv-riscalate} and the Lebesgue's theorem, we have,
for every $x \in [0,1]$,
\begin{equation} \label{n7}
\lim_{\ep \to 0^+} F_{\ep ,L}(x)=\dfrac{a'(x)}{a(x)^2}\dfrac{1}{L} \int_{-L/2}^{L/2}\mathcal{L}(\bar{v}' (z))\,dz,
\end{equation}
where $\bar{v}$ is the limit
of $v_{x,\ep}$ given by Proposition \ref{conv-riscalate}.
\vs{4}
\ni
If $\bar{v}$ is a constant function (a situation which can occur only
if $E(x) = 0$ or $E(x) = W_0$), we have already concluded. Indeed, in this case the above integral
equals zero, and $K(E(x)) = 0$ as well.
Otherwise, we know that
$\bar{v}(z) = \bar{v}_{a(x),E(x)}(\pm z + t(x))$ for some $t(x) \in \mathbb{R}$ and we can write
$$
\begin{array}{l}
\displaystyle{\inf_{t\in \RR}\dfrac{a'(x)}{a(x)^2} \dfrac{1}{L} \left(\int_{-L/2}^{L/2}\mathcal{L}(\bar{v}'_{a(x),E(x)}(\pm z+t))\,dz \right)\psi(x)\leq \liminf_{\ep \to 0^+} F_{\ep ,L}(x)\psi (x)}\\
\\
\displaystyle{\leq \limsup_{\ep \to 0^+} F_{\ep ,L}(x)\psi (x)\leq \int_0^1\sup_{t\in \RR} \dfrac{a'(x)}{a(x)^2}\dfrac{1}{L} \left(\int_{-L/2}^{L/2}\mathcal{L}(\bar{v}'_{a(x),E(x)}(\pm z+t))\,dz \right)\psi(x)}.
\end{array}
$$
Using Fatou's Lemma, we thus get
$$
\begin{array}{l}
\displaystyle{\int_0^1\inf_{t\in \RR}\dfrac{a'(x)}{a(x)^2} \dfrac{1}{L} \left(\int_{-L/2}^{L/2}\mathcal{L}(\bar{v}'_{a(x),E(x)}(\pm z+t))\,dz \right)\psi(x)\,dx\leq \liminf_{\ep \to 0^+} \int_0^1 F_{\ep ,L}(x)\psi (x)\,dx}\\
\\
\displaystyle{\leq \limsup_{\ep \to 0^+}\int_0^1F_{\ep ,L}(x)\psi (x)\,dx\leq \int_0^1 \sup_{t\in \RR} \dfrac{a'(x)}{a(x)^2}\dfrac{1}{L} \left(\int_{-L/2}^{L/2}\mathcal{L}(\bar{v}'_{a(x),E(x)}(\pm z+t))\,dz \right)\psi(x)\,dx}.
\end{array}
$$
Accordingly, the proof can be concluded by showing that,
\begin{equation} \label{n11}
\lim_{L\to +\infty} \dfrac{1}{L} \int_{-L/2}^{L/2} \mathcal{L}(\bar{v}'_{a(x),E(x)}(\pm z+t))\,dz =
a(x)K(E(x)),
\end{equation}
uniformly in $x \in [0,1]$ and $t\in \RR$.
\vs{4}
\ni
To see this, we distinguish two cases. If $E(x) = 0$, then $K(E(x)) = 0$ and the left-hand side is zero as well,
since 
$$
\lim_{L \to +\infty} \int_{-L/2}^{L/2} \mathcal{L}(\bar{v}'_{a(x),E(x)}(\pm z+t))\,dz = 
\int_{-\infty}^{+\infty}  \mathcal{L}(\bar{v}'_{a(x),E(x)}(\pm z))\,dz < +\infty,
$$
in view of Proposition \ref{eterocline}. If $E(x) > 0$, we 
set for simplicity of notation $T_x = T_{a(x)}(E(x))$ and 
write $L=n_L T_x+r_L$, with $n_L\in \NN$ and $r_L\in [0,T_x)$; then, we have
\[
\begin{array}{ll}
&\displaystyle{\dfrac{1}{L} \int_{-L/2}^{L/2} \mathcal{L}(\bar{v}'_{a(x),E(x)}(\pm z+t))\,dz}
\\= &\displaystyle{\dfrac{1}{n_L T_x +r_L} \int_{-n_L T_x/2}^{n_L T_x/2}\mathcal{L}(\bar{v}'_{a(x),E(x)}(\pm z+t))\,dz }
\vspace{0.1cm}
\\ & +\displaystyle{\dfrac{1}{n_L T_x +r_L} \int_{-(n_L T_x+r_L)/2}^{-n_L T_x/2}\mathcal{L}(\bar{v}'_{a(x),E(x)}(\pm z+t))\,dz}
\vspace{0.1cm}
\\
& +\displaystyle{\dfrac{1}{(n_L T_x +r_L)/2} \int_{n_L T_x/2}^{(n_L T_x+r_L)/2}\mathcal{L}(\bar{v}'_{a(x),E(x)}(\pm z+t))\,dz}\\
\\= &\displaystyle{\dfrac{n_L}{n_L T_x +r_L} \int_{0}^{T_x}\mathcal{L}(\bar{v}'_{a(x),E(x)}(z))\,dz + \ldots.}
\end{array}
\]
Now, the first integral above has limit, for $L \to +\infty$, equal to 
$K_{a(x)}(E(x))$, while the remainder is easily seen to go to zero.
We can thus conclude using Lemma \ref{relaz-periodi-kappa}.
\end{proof}
\vs{12}
\ni
We are now interested in the existence of solutions of
\eqref{eq-profilo} vanishing somewhere on $[0,1]$
(notice that this is possible since - from Proposition \ref{prop-kappa} -
$K$ is not Lipschitz-continuous at $\xi = 0$).
More precisely, we have the following proposition, which is proved
in \cite[Proposition 1.3]{FeMaTa-05} (notice, indeed, that the properties of $K$ collected in Proposition
\ref{prop-kappa} remain the same with respect to the case $p = 2$).
\vs{12}
\ni
\begin{proposition}\label{prop-soluzioni}
Let $E: [0,1] \to [0,W_0]$ be a solution of \eqref{eq-profilo}, vanishing somewhere on $[0,1]$.
Then, the connected components of $\{ x : E(x) > 0\}$ are intervals of the following type:
\begin{itemize}
\item[(i)] $(s,t)$, where $0 \leq s < t \leq 1$ satisfy $a(s) = a(t)$ and $a(x) > a(s)$ for $x \in (s,t)$,
\item[(ii)] $(s,1]$, where $s \in [0,1)$ satisfies $a(x) > a(s)$ for $x \in (s,1]$,
\item[(iii)] $[0,t)$, where $t \in (0,1]$ satisfies $a(x) > a(t)$ for $x \in [0,t)$.
\end{itemize}
Conversely, if $A \subset [0,1]$ is a disjoint union of intervals of the type
(i), (ii), (iii), there exists a solution $E$ of \eqref{eq-profilo} such that
$E(x) = 0$ if and only if $x \in [0,1]\setminus A$. Moreover, if $A \neq [0,1]$
such a solution is unique.
\end{proposition} 
\vs{12}
\ni
\section{On the distribution of zeros}\label{section4}
\def\theequation{4.\arabic{equation}}\makeatother
\setcounter{equation}{0}
\vs{12}
\ni
In this section we prove some results about the asymptotic distribution of zeros of $u_\ep$ in $[0,1]$,
when $\ep \to 0^+$. Both in Propositions \ref{stima-zeri} and \ref{accumulazione-zeri} below, we suppose that
$\ep_n \to 0^+$ is a sequence such that the
energy $E_n = E_{\ep_n}$ of $u_n = u_{\ep_n}$ converges to some 
limit $E$.
\vs{12}
\ni
\begin{proposition} \label{stima-zeri} 
Let us denote by $z_{n}$ the number of zeros of $u_n$ in $[0,1]$. 
Then, the following relation holds true:
\begin{equation} \label{zeri-stima}
\lim_{n \to +\infty} \ep_n z_n =\int_0^1 \dfrac{2\sqrt[p]{a(x)}}{T(E(x))}\,dx,
\end{equation}
where the right-hand side of \eqref{zeri-stima} has to be considered equal to zero when $E(x)\equiv 0$.
\end{proposition}
\vs{8}
\ni
\begin{proof}
Let us consider the function $g:[-1,1]\to \RR$ defined by
\[
g(u)=\sqrt[p]{p}\left(W_0-W(u)\right)^{1/p}\, \textnormal{sgn}(u),\quad \forall \ u\in [-1,1].
\]
Due to the assumption (W1), $g$ is of class $C^1$ and the differential equation in \eqref{bvp-1} can be written as
\begin{equation} \label{eq-riscritta}
\ep (\phi_p(\ep u'))'+a(x)\phi_p(g(u))g'(u)=0.
\end{equation}
Let us make the change of variables
\begin{equation} \label{cambio-var-2}
\left\{\begin{array}{l}
       a^{1/p} g(u)=r^{2/p} C_p(\theta)\\
		   \\
			\phi_p(\ep u')=r^{2/p^*} S_p(\theta),
			\end{array}
\right.
\end{equation}
where $C_p$ and $S_p$ are, respectively, the \emph{p-cosine} and
the \emph{p-sine} functions (see \cite{DeElMa-89}). We recall that such functions satisfy the following properties:
\begin{itemize}
\item[(i)] $C_p$ and $S_p$ are $2\pi_p$-periodic,
\item[(ii)] $C_p(\theta) = 0$ if and only if $\theta = \pi_p/2 + k \pi_p$ for some $k \in \mathbb{Z}$,
\item[(iii)] $C_p'(\theta) = - \phi_{p^*}(S_p(\theta))$ and $S_p'(\theta) = \phi_p(C_p(\theta))$,
\item[(iv)] $\vert C_p(\theta) \vert^p /p + \vert S_p(\theta) \vert^{p^*}/p^* \equiv 1/p$.
\end{itemize}
The change of variable in \eqref{cambio-var-2} is admissible, since
it is well-known (see \cite[Lemma 2.1]{MaNjZa-95}) that nontrivial solutions of
\eqref{eq-riscritta} have only simple zeros (i.e., $u(x)^2 + u'(x)^2 \neq 0$); moreover, due 
to (ii), if $z$ denotes the number of (simple) zeros of a nontrivial solution of
\eqref{eq-riscritta}, then we have
\begin{equation}\label{zerirotazioni}
\left|z - \dfrac{\theta (1)-\theta(0)}{\pi_p}\right|\leq 1;
\end{equation}
To get an estimate of the above quantity, we argue as follows.
By differentiating the first relation in \eqref{cambio-var-2}, we obtain
$$
\dfrac{\ep}{p} a^{1/p-1}a'g(u)+\ep a^{1/p}g'(u)u'=\dfrac{2\ep}{p}r^{2/p-1}r'C_p(\theta)+\ep r^{2/p}C_p'(\theta)\theta';
$$
that is - using (iii) and the fact that $\ep u' = \phi_{p^*}\left( r^{2/p^*} S_p(\theta)\right)$ - 
\begin{equation} \label{b11}
\dfrac{2\ep }{p}r^{2/p-1}r'C_p(\theta)-\ep r^{2/p}\phi_{p^*}(S_p(\theta))\theta'=\dfrac{\ep}{p} a^{1/p-1}a'g(u)+a^{1/p}g'(u)
\phi_{p^*}\left(r^{2/p^*}S_p(\theta)\right).
\end{equation}
On the other hand, differentiating the second relation in \eqref{cambio-var-2} and using
the equation \eqref{eq-riscritta}, we infer
$$
\dfrac{2\ep}{p^*}r^{2/p^*-1}r'S_p(\theta)+\ep r^{2/p^*}S_p'(\theta)\theta'=-a \phi_p(g(u))g'(u),
$$
that is - using (iii) and the fact that $a^{p^*} \phi_p(g(u)) = \phi_p \left( r^{2/p}C_p(\theta)\right)$ -
\begin{equation} \label{b12}
\dfrac{2\ep }{p^*}r^{2/p^*-1}r'S_p(\theta)+\ep r^{2/p^*}\phi_{p}\left(C_p(\theta)\right)\theta'=-a^{1/p}g'(u)\phi_{p}\left(r^{2/p}C_p(\theta)\right).
\end{equation}
We now multiply \eqref{b12} by $r^{2/p}C_p(\theta)/p$ and \eqref{b11} by ${r^{2/p^*}S_p(\theta)/p^*}$ and subtract them,
obtaining - in view of (iv) - 
$$
\ep r^2 \theta' = - r^2 a^{1/p} g'(u) - \frac{\ep}{p p^*} r^{2/p^*} S_p(\theta) a' g(u),
$$
which can be rewritten - using the fact that $g(u)/r^{2/p} = C_p(\theta)/a^{1/p}$ - as
$$
\ep \theta' = - a^{1/p} g'(u) - \frac{\ep}{p p^*} \frac{a'}{a^{1/p}} C_p(\theta) S_p(\theta).
$$
\vs{4}
\ni
We can now pass to prove \eqref{zeri-stima}. 
In view of \eqref{zerirotazioni}, we have
\begin{equation} \label{b14}
\lim_{n \to +\infty} \ep_n z_n=-\lim_{n \to +\infty} \dfrac{1}{\pi_p}\int_0^1 \ep_n \theta_n'(x)\,dx=
\lim_{n \to +\infty} \dfrac{1}{\pi_p}\int_0^1 \sqrt[p]{a(x)}\, g'(u_n (x))\,dx.
\end{equation}
Now, for every $x\in (0,1)$ let us consider the function $v_{x,n}$ defined as in \eqref{def-riscalate0},
that is, $v_{x,n}(y) = u_n(x+ \ep_n y)$, and, according to Proposition \ref{conv-riscalate},
let $\bar{v}$ be the limit of $v_{x,n}$ for $n \to +\infty$.
\vs{4}
\ni
If $\bar{v}$ is a constant function, we easily conclude. Indeed, if $\bar{v} = 0$ (hence
$E(x) = W_0$), then
\[
\lim_{n \to +\infty} \int_0^1 \sqrt[p]{a(x)}g'(u_n(x))\,dx=g'(0)\ \int_0^1 \sqrt[p]{a(x)}\,dx=\sqrt[p]{C_0} \ \int_0^1 \sqrt[p]{a(x)}\,dx
\]
and - from Proposition \ref{prop-periodo} - 
\[
\lim_{\xi \to W_0^-} T(\xi)=\dfrac{2\pi_p}{\sqrt[p]{C_0}}.
\]
On the other hand, if $\bar{v} = \pm 1$ (hence $E(x) = 0$), then
\[
\lim_{n \to +\infty} \int_0^1 \sqrt[p]{a(x)}g'(u_n(x))\,dx=g'(0)\ \int_0^1 \sqrt[p]{a(x)}\,dx=0
\]
and - from Proposition \ref{prop-periodo} - 
\begin{equation}\label{timemapaux}
\lim_{\xi \to 0^+} T(\xi)=+\infty.
\end{equation}
\vs{4}
\ni
If $\bar{v}$ is non-constant, we have to argue similarly as in the proof of Theorem 
\ref{profilo-energia}, using mollifiers and convolution.
In particular, if $\bar{v}$ is the heteroclinic (hence, $E(x) = 0$) we can prove that
\[
\lim_{n \to +\infty} \int_0^1 \sqrt[p]{a(x)}g'(u_n(x))\,dx=0
\]
and we conclude again in view of \eqref{timemapaux}. On the other hand, if 
$\bar{v}$ is periodic (hence, $0 < E(x) < W_0$) with minimal period
$T_{a(x)}(E(x))$, we have 
\begin{equation} \label{b15}
\lim_{\ep \to 0^+} \dfrac{1}{\pi_p}\int_0^1 \sqrt[p]{a(x)}\, g'(u_\ep (x))\,dx=\int_0^1 \dfrac{1}{\pi_p} \sqrt[p]{a(x)}\dfrac{1}{T_x} \int_0^{T_x} g'(\bar{v}(y))\,dy\,dx.
\end{equation}
By writing the conservation of energy for $\bar{v}$ as
$$
\mathcal{L}(\bar{v}') + \frac{a(x)}{p}\vert g(\bar{v}) \vert^p = a(x)(W_0-E(x)),
$$
and changing variables via $u = g(\bar{v}(y))$ in \eqref{b15}, 
we obtain
$$
\dfrac{1}{T_x} \int_0^{T_x} g'(\bar{v}(y))\,dy=\dfrac{2}{T_x} \int_{k_-(E(x))}^{k_+(E(x))} \dfrac{du}{\LL_+^{-1} \left(a(x)\left(W_0-E(x)-|u|^p/p\right)\right)},
$$
where $k_{\pm}(s)=\pm \sqrt[p]{p(W_0-s)}$, for every $s\in (0,W_0)$.
Using \eqref{relaz-periodi-kappa} and elementary computations
\begin{align*}
\dfrac{1}{T_x} \int_0^{T_x} g'(\bar{v}(y))\,dy & = \dfrac{2}{T_x \,a(x)^{1/p }(p^*)^{1/p}} 
\int_{k_-(E(x))}^{k_+(E(x))} \dfrac{du}{\left(W_0-E(x)-|u|^p/p\right)^{1/p}}\\
& =\dfrac{2}{T(E(x))} \dfrac{(W_0-E(x))^{-1/p}}{(p^*)^{1/p}}\int_{k_-(E(x))}^{k_+(E(x))} \dfrac{du}{\left(1-\left(\dfrac{|u|}{\sqrt[p]{p(W_0-E(x))}}\right)^p\right)^{1/p}}
\\
& =\dfrac{2}{T(E(x))} \left(\dfrac{p}{p^*}\right)^{1/p}\int_{-1}^{1} \dfrac{ds}{\left(1-|s|^p\right)^{1/p}}
=\dfrac{2\pi_p}{T(E(x))}.
\end{align*}
Recalling \eqref{b15}, this concludes the proof.
\end{proof}
\vs{12}
\ni
We conclude this section with a result describing the set of accumulations of zeros of
$u_n$, in connection with the support of the limit function $E$ and the set of critical points
of the weight $a$. 
\vs{12}
\ni
\begin{proposition}\label{accumulazione-zeri}
Let us denote by $Z$ the set of accumulations of zeros of $u_{n}$, that is,
$$
Z = \cap_{n = 1}^{\infty} \cup_{j=n}^{\infty} \overline{\left\{ x \in [0,1] : u_{j}(x) = 0\right\}}.
$$ 
Then, the following inclusions hold true:
$$
\textnormal{supp}\, E \subset Z \subset \textnormal{supp}\, E \cup \{ x : a'(x) = 0\} \cup \{0,1\}.
$$
\end{proposition}
\vs{12}
\ni
This result follows from the next proposition, which will also be used in the proof of Theorem \ref{main-ex-1}.
\vs{12}
\ni
\begin{proposition}\label{prop2.6}
Let $u_n$ be a family of solutions of \eqref{bvp-1} such that 
\[
E_n\to 0,\quad \mbox{uniformly in $[0,1]$}.
\]
Suppose that $[\alpha,\beta]\subset [0,1)$ satisfies, for some $h>0$, 
\[
a'(x)>0,\quad \forall \ x\in [\max(0,\alpha-h),\beta+h],
\]
and
\[
E(x)=0,\quad \forall \ x\in [\max(0,\alpha-h),\beta+h].
\]
Then, for $n$ sufficiently large it holds
\[
u_n'(x)\neq 0,\quad \forall \ x\in (\alpha,\beta].
\]
A symmetric result holds true in the case $[\alpha,\beta] \subset (0,1]$.
\end{proposition}
\vs{12}
\ni
Proposition \ref{prop2.6} follows from the variational characterization of solutions developed
in Section \ref{section2}. Since the complete proof is very long, but
requires only minor modifications with respect to the case $p=2$
treated in \cite[Proposition 2.6]{FeMaTa-05},
we omit it.
\vs{12}
\ni
\section{Existence of highly oscillatory solutions}\label{section5}
\def\theequation{5.\arabic{equation}}\makeatother
\setcounter{equation}{0}
\vs{12}
\ni
In this section we prove the existence of solutions $u_\ep$ of \eqref{bvp-1} such that 
\[
E_\ep \to E,\quad \mbox{uniformly in $[0,1]$},
\]
for a given energy profile $E$ satisfying \eqref{eq-profilo}.
\vs{6}
\ni
To this aim we use a broken-geodesic approach. We thus consider, for $\ep > 0$ and $[s,t] \subset [0,1]$, the energy functional
$$
I_\ep(s,t;u) = \int_s^t \left(\frac{\ep^p}{p} \vert u'(x) \vert^p + a(x)W(u(x))\right)\,dx, \qquad u \in W^{1,p}(s,t),
$$
assuming that $W$ is extended outside $[-1,1]$ as
\begin{equation} \label{estensione}
\begin{array}{l}
\displaystyle{W(u)=\dfrac{C_1}{p} (u-1)^p,\quad u\geq 1}\\
\\
\displaystyle{W(u)=\dfrac{C_{-1}}{p} |u+1|^p,\quad u\leq -1.}
\end{array}
\end{equation}
Let us observe that the extended function $W$ satisfies $(W2)$ in $\RR$.
\vs{4}
\ni
For every $[s,t]\subset (0,1)$ we define
\begin{equation}\label{mpiudir}
m_+(\ep;s,t) = \inf \left\{ I_\ep(s,t;u) \, : \, u \in W^{1,p}(s,t), \, u \geq 0, \, u(s) = u(t) = 0 \right\}
\end{equation}
and
\begin{equation}\label{mmenodir}
m_-(\ep;s,t) = \inf \left\{ I_\ep(s,t;u) \, : \, u \in W^{1,p}(s,t), \, u \leq 0, \, u(s) = u(t) = 0 \right\}.
\end{equation}
Analogously, for $[s,1]\subset (0,1]$ and $[0,t]\subset [0,1)$, let
\begin{equation}\label{mpiuzero}
m_+(\ep;0,t) = \inf \left\{ I_\ep(0,t;u) \, : \, u \in W^{1,p}(0,t), \, u \geq 0, \, u(t) = 0 \right\}
\end{equation}
and
\begin{equation}\label{mpiuuno}
m_+(\ep;s,1) = \inf \left\{ I_\ep(0,t;u) \, : \, u \in W^{1,p}(s,1), \, u \geq 0, \, u(s) = 0 \right\}.
\end{equation}
The critical levels $m_-(\ep;0,t)$ and $m_-(\ep;s,1)$ are defined replacing $u\geq 0$ with $u\leq 0$.

\ni
We will show (see Lemma \ref{minimidir} and Lemma \ref{minimidir2}) that all these minimization problems have a unique minimizer $u_{\pm}$.
\vs{12}
\ni
We now state a first result on the existence of solutions of \eqref{bvp-1}, corresponding to the case in which the 
the set $\{ x : E(x) > 0\}$ is the union of two disjoint intervals $(s_0,t_0)$ and $(s_1,t_1)$ 
(see Proposition \ref{prop-soluzioni}) satisfying 
\begin{equation} \label{num2}
a'(s_i)>0,\qquad a'(t_i)<0,\quad i=0, 1.
\end{equation}
Let $h_0>0$ such that $a'>0$ in $[s_0-h_0,s_0]\cup [s_1-h_0,s_1]$ and $a'<0$ in $[t_0,t_0+h_0]\cup [t_1,t_1+h_0]$ and let $n^i_\ep$, $i=0, 1$, be positive integers such that
\begin{equation} \label{ipotesinumerozeri}
\ep n^i_\ep \to \int_{s_i}^{t_i}\dfrac{2\sqrt[p]{a(x)}}{T(E(x))}\,dx,\quad \ep\to 0^+.
\end{equation}
Let us also consider
\[
\begin{array}{ll}
\Delta'=&\{(\tau_1,\ldots,\tau_{n^0_\ep+n^1_\ep}):\ s_0-h_0\leq 
\tau_1\leq \tau_2\leq \ldots\leq \tau_{n^0_\ep}\leq t_0+h_0,\\
&\\
&s_1-h_0\leq \tau_{n^0_\ep+1}\leq  \ldots\leq \tau_{n^0_\ep+n^1_\ep}\leq t_1+h_0\}; 
\end{array}
\]
for every $\ep>0$ we define
\[
f_\ep(\tau_1,\ldots,\tau_{n^0_\ep+n^1_\ep})=\sum_{j=0}^{n^0_\ep+n^1_\ep} m_{(-)^j}(\ep;\tau_j,\tau_{j+1}),\quad \forall  \ (\tau_1,\ldots,\tau_{n^0_\ep+n^1_\ep})\in \Delta',
\]
where for every $k\in \NN$ we have
\[
(-)^{2k}=+,\quad (-)^{2k-1}=-
\]
and $\tau_0=0$, $\tau_{n^0_\ep+n^1_\ep+1}=1$.
\vs{6}
\ni
We then have the following result:
\vs{12}
\ni
\begin{theorem} \label{main-ex-1} For every $\ep$ sufficiently small, the maximization problem
\[
\max_{(\tau_1,\ldots,\tau_{n^0_\ep+n^1_\ep})\in \Delta'} f_\ep (\tau_1,\ldots \tau_{n^0_\ep+n^1_\ep})
\]
has a maximizer $(\tau_1,\ldots, \tau_{n^0_\ep+n^1_\ep})\in \mathring{\Delta}'$ such that the corresponding minimizer $u_{(-)^j}(\ep;\tau_j,\tau_{j+1})$ of $m_{(-)^j}(\ep;\tau_j,\tau_{j+1})$ is nontrivial, for every $j=0,\ldots,n^0_\ep+n^1_\ep$. 
\vs{4}
\ni
Moreover, the function $u_\ep:[0,1]\to \RR$ defined by
\begin{equation} \label{num}
u_\ep (x)=u_{(-)^j}(\ep;\tau_j,\tau_{j+1})(x),\quad x\in [\tau_j,\tau_{j+1}],
\end{equation}
is a solution of \eqref{bvp-1} such that
\[
E_\ep \to E,\quad \mbox{uniformly in $[0,1]$}.
\]
Finally
\[
\ep n_\ep ([s_i,t_i])\to \int_{s_i}^{t_i}\dfrac{2\sqrt[p]{a(x)}}{T(E(x))}\,dx,\quad \ep\to 0^+,\ i=0, 1,
\]
where $n_\ep ([s_i,t_i])$ is the number of zeros of $u_\ep$ in $[s_i,t_i]$, $i=0, 1$.
\end{theorem}
\vs{12}
\ni
In order to prove Theorem \ref{main-ex-1} we need several preliminary results
dealing with the minimization problems \eqref{mpiudir}, \eqref{mpiuzero}, \eqref{mpiuuno}.
We state them in the case of $m_+$ and $u_+$, but analogous conclusions hold true for 
$m_-$ and $u_-$. 
\vs{12}
\ni
\begin{lemma}\label{minimidir}
The following results hold true:
\begin{itemize}
\item[(i)] the minimization problem \eqref{mpiudir} has a unique minimizer $u_+(\ep;s,t)$;
\item[(ii)] $u'_+(\ep;s,t)(s)$, $u'_+(\ep;s,t)(t): \{ (s,t): \, 0 < s < t < 1 \} \to \mathbb{R}$ are continuous;
\item[(iii)] $m_+(\ep;s,t)$ is differentiable with respect to $s$ and $t$
in $\{ (s,t): \, 0 < s < t < 1 \}$ and
$$
\frac{\partial }{\partial s} m_+(\ep;s,t) = \frac{\ep^p}{p^*}\vert u'_+(\ep;s,t)(s) \vert^p - a(s) W_0
$$
$$
\frac{\partial }{\partial t} m_+(\ep;s,t) = - \frac{\ep^p}{p^*}\vert u'_+(\ep;s,t)(t) \vert^p + a(t) W_0.
$$
\end{itemize}
\end{lemma}
\vs{8}
\ni
\begin{proof}
(i) The existence of a minimizer is straightforward. Indeed, 
the boundary condition $u(t) = 0$ implies that, for any $x \in [s,t]$,
$$
\vert u(x) \vert \leq \left\vert \int_x^t \vert u' \vert \,dx\right\vert
\leq \vert t - x \vert^{1/p^*} \left\vert \int_x^t \vert u' \vert^p \,dx\right\vert^{1/p}.
$$
Hence, $I_\ep(s,t;u)$ is coercive (and weakly lower semicontinuous)
on a convex subset of the reflexive Banach space $W^{1,p}(s,t)$, and the direct method of the Calculus of Variations
applies.

\ni
As far as the uniqueness is considered, it is sufficient to apply the results of \cite{BrOs-86} and \cite{BeKa-02}; indeed,
by $(W2)$, the function $W'(u)/\phi_p(u)$ is increasing in $(0,1]$.
\vs{6}
\ni
(ii-iii) These result are quite standard; they are based on the uniqueness of minimizers proved in (i)
together with simple calculations.
\end{proof}
\vs{12}
\ni
In a similar way it is possible to prove the following result.
\vs{12}
\ni
\begin{lemma} \label{minimidir2}
The following results hold true:
\begin{itemize}
\item[(i)] the minimization problems \eqref{mpiuzero}, \eqref{mpiuuno} have unique minimizers
$u_+(\ep;0,t)$ and $u_+(\ep;s,1)$;
\item[(ii)] $u_+(\ep;0,t)(0)$, $u_+'(\ep;0,t)(t)$, $u_+(\ep;s,1)(1)$, $u_+'(\ep;s,1)(s)$ are continuous functions;
\item[(iii)] $m_+(\ep;0,t)$ and $m_+(\ep;s,1)$ are differentiable with respect to $t$ and $s$, respectively, 
and
$$
\frac{\partial }{\partial t} m_+(\ep;0,t)= - \frac{\ep^p}{p^*}\vert u'_+(\ep;0,t)(t) \vert^p + a(t) W_0
$$
$$
\frac{\partial }{\partial s} m_+(\ep;s,1) = \frac{\ep^p}{p^*}\vert u'_+(\ep;s,1)(s) \vert^p - a(s) W_0.
$$
\end{itemize}
\end{lemma}
\vs{12}
\ni
The next result gives some knowledge on the derivatives of the critical levels $m_\pm$
in connection with the monotonicity of the weight function $a$.
\vs{12}
\ni
\begin{lemma} \label{stimeminimi}
Suppose that $a'(x) > 0$ for every $x \in [\alpha,\beta]$. There exists
$C_0 > 0$ such that, for $\ep > 0$ small enough,
\begin{itemize}
\item[(i)] for $\alpha \leq s < t \leq \beta$,
$$
\frac{\partial}{\partial t} m_+(\ep;s,t) > 0;
$$
\item[(ii)] for $\alpha \leq s < t \leq \beta$,
$$
\left( \frac{\partial}{\partial s} + \frac{\partial}{\partial t} \right) m_+(\ep;s,t) > 0;
$$
\item[(iii)] for $s \in [\alpha,\beta]$ and $t \in (s,1]$
$$
\frac{\partial}{\partial s} m_+(\ep;s,t) > 0, \qquad \mbox{if } \quad\frac{t-s}{\ep} \geq C_0 \vert \log \ep \vert;
$$
for $t \in [\alpha,\beta]$ and $s \in (0,t]$
$$
\frac{\partial}{\partial t} m_+(\ep;s,t) > 0, \qquad \mbox{if } \quad \frac{t-s}{\ep} \geq C_0 \vert \log \ep \vert.
$$
\end{itemize}
\end{lemma}
\vs{8}
\ni
\begin{proof}
We give the proof when $0 < \alpha < \beta < 1$, the other cases being similar.
As a preliminary observation, we notice that we have  
\begin{equation}\label{m-con-energia}
\frac{\partial}{\partial t} m_+(\ep;s,t) = a(t)E(t), \qquad \frac{\partial}{\partial s} m_+(\ep;s,t) = -a(s)E(s),
\end{equation}
where $E(x)$ is defined, as in \eqref{def-energia}, 
to be the energy of the function $u_+(\ep;s,t)$.
\vs{6}
\ni
(i) Let $x_0 \in (s,t)$ be the maximum point of $u_+(\ep;s,t)$; then,
$E(x_0) = W(u_p(\ep;s,t)(x_0)) > 0$. Using \eqref{eq-energia}, and since
$a'(x) > 0$ for $x \in [s,t]$, we have $E(t) > E(x_0)$.
Recalling \eqref{m-con-energia}, this concludes the proof.
\vs{6}
\ni
(ii) Again from \eqref{m-con-energia} and \eqref{eq-energia},
\begin{align*}
\left( \frac{\partial}{\partial s} + \frac{\partial}{\partial t} \right) m_+(\ep;s,t) & = \int_s^t \frac{d}{dx}(a(x)E(x))\,dx \\
& = \int_s^t a'(x)W(u_+(\ep;s,t)(x))\,dx > 0.
\end{align*} 
\vs{6}
\ni
(iii) We give the proof for $\tfrac{\partial}{\partial s} m_+(\ep;s,t)$.
Define, for $y \in [0,\tfrac{t-s}{\ep}]$,
$$
v_\ep(y) = u_+(\ep;s,t)(s+\ep y);
$$
moreover, set $\mu_\ep = \tfrac{1}{K_2} \vert\log\ep \vert$, with $K_2 > 0$
the constant appearing in \eqref{stima-asintotica-Nakashima-2}. Then, 
provided
$$
\frac{t-s}{\ep} \geq 2\mu_e,
$$
it holds that
\begin{equation}\label{a8}
\left( 1 - v_\ep(y) \right) + \left\vert v_\ep'(y) \right\vert \leq C_1 \ep, \quad \mbox{ for }
\; y \in \left[ \mu_\ep,\frac{t-s}{2\ep}\right].
\end{equation}
Now, let $\varphi \in C^1(\mathbb{R}^+)$ be a decreasing function such that
$$
\varphi(z) = 1 \quad\forall z \in [0,1], \qquad \varphi(z) = 0 \quad \forall z \geq 2.
$$
We have, using \eqref{m-con-energia} and \eqref{eq-energia},
\begin{align*}
\frac{1}{\ep}\frac{\partial}{\partial s} m_+(\ep;s,t) & = \frac{1}{\ep}
\int_s^t \frac{d}{dx}\left[ \varphi\left( \frac{x-s}{\ep\mu_\ep}\right)a(x)E(x)\right]\,dx \\
& = \frac{1}{\ep^2}\int_s^t \frac{1}{\mu_\ep}\varphi'\left( \frac{x-s}{\ep\mu_\ep}\right)a(x)E(x)\,dx 
\\ & + \frac{1}{\ep}\int_s^t \varphi\left( \frac{x-s}{\ep\mu_\ep}\right)a'(x)W(u_+(\ep;s,t)(x))\,dx 
\\ & = \int_{\mu_\ep}^{2\mu_\ep} \frac{1}{\ep\mu_\ep}\varphi'\left( \frac{y}{\mu_\ep}\right)a(s+\ep y)E(s+ \ep y)\,dy 
\\ & + \int_0^{2\mu_\ep} \varphi\left( \frac{y}{\mu_\ep}\right)a'(s+\ep y)W(v_\ep(y))\,dy. 
\end{align*}
Now, observe that from assumption $(W1)$ it follows that, for a suitable $C_3 > 0$,
$$
\vert W(u) \vert \leq C_3 \vert 1 - u \vert^p, \quad \mbox{ for every } u \in \mathbb{R}.
$$
As a consequence, \eqref{a8} implies that $\vert E(s + \ep y) \vert \leq C_4 \ep^p$, so that
$$
\left\vert \int_{\mu_\ep}^{2\mu_\ep} \frac{1}{\ep\mu_\ep}\varphi'\left( \frac{y}{\mu_\ep}\right)a(s+\ep y)E(s+ \ep y)\,dy \right \vert \leq C_4 \ep^{p-1}.
$$
On the other hand, recalling Proposition \ref{conv-riscalate}, we have that, up to subsequences, $v_\ep$ converges locally uniformly to a solution ${\bar v}$ of 
\[
-(\phi_p({\bar v}'))'+a(s)W({\bar v})=0
\]
such that ${\bar v}(0)=0$; hence we deduce that there exists $\delta=\delta(s)>0$ such that
\begin{equation} \label{nnn1}
W({\bar v}(y))\geq \delta,\quad \forall \ y\in [0,1].
\end{equation}
Since $a$ is striclty positive in $[0,1]$, it is possible to choose $\delta$ indipendent on $s$ in \eqref{nnn1}.

\ni
As a consequence, there exists $\delta'>0$ such that
\[
\left\vert \int_0^{2\mu_\ep} \varphi\left( \frac{y}{\mu_\ep}\right)a'(s+\ep y)W(v_\ep(y))\,dy \right \vert\geq \delta'>0,
\]
for every $\ep>0$ sufficiently small.
This concludes the proof.
\end{proof}
\vs{12}
\ni
Our last results specify when the minimizers $u_\pm$ are non trivial.
\vs{12}
\ni
\begin{lemma} \label{non-banali}
For $[s,t] \subset (0,1)$, $m_+(\ep;s,t) < W_0\int_s^t a(x)\,dx$ if and only if 
$$
\inf_{u \in W^{1,p}_0(s,t)} \frac{\int_s^t \left( \ep^p\vert u' \vert^p -  C_0 a(x) \vert u \vert^p\right)}
{\int_s^t \vert u \vert^p} < 0,
$$
i.e., the first eigenvalue of the problem
\begin{equation}\label{princeigen}
-\ep(\phi_p(\ep u'))'- C_0 a(x)\phi_p(u)= \lambda \phi_p(u), \qquad u \in W^{1,p}_0(s,t)
\end{equation}
is negative. In this case, $u_+(\ep;s,t)(x) > 0$ for every $x \in (s,t)$.
\end{lemma}
\vs{8}
\ni
\begin{proof}
Using the assumption on $W$ we have 
\begin{align*}
W(u) - W_0 & = \int_0^u W'(v)\,dv \geq \int_0^u \frac{W'(v)}{\phi_p(v)}{\phi_p(v)}\,dv 
\\ & \geq - C_0\int_0^u \phi_p(v)\,dv = - \dfrac{C_0}{p} \vert u \vert^p,
\end{align*}
so that
$$
I_\ep(s,t;u) - W_0 \int_s^t a(x)\,dx \geq \int_s^t \left( \frac{\ep^p}{p}\vert u' \vert^p - \dfrac{C_0}{p} a(x) \vert u \vert^p\right).
$$
Hence, if $m_+(\ep;s,t) < W_0\int_s^t a(x)\,dx$ then $u = u_+(\ep;s,t)$ satisfies
$$
\int_s^t  \left( \ep^p\vert u' \vert^p -  C_0 a(x)\vert u \vert^p\right) < 0.
$$
Conversely, let $e(x)$ be the first (positive) eigenfunction of \eqref{princeigen};
then
\begin{align*}
\lim_{h \to 0}\frac{I_\ep(s,t;\pm h e) - I_\ep(s,t;0)}{\vert h \vert^p} & = 
\lim_{h \to 0}\int_s^t \left(  \frac{\ep^p}{p}\vert e' \vert^p + a(x)\frac{W(he) - W_0}{\vert he \vert^p}\vert e \vert^p\right) \\
& = \int_s^t \left(  \frac{\ep^p}{p}\vert e' \vert^p - \dfrac{C_0}{p}
 a(x)\vert e \vert^p\right) < 0.
\end{align*}
Hence, if $\vert h \vert$ is small enough, $I_\ep(s,t;\pm he) < I_\ep(s,t;0) = W_0 \int_s^t a(x)\,dx$.
\\
\ni
The fact that $u_+ > 0$ is easily checked.
\end{proof}
\vs{12}
\ni
\begin{lemma} \label{A6} Let $l_\ep(s,t)$ be the number of negative eigenvalues of
\begin{equation} \label{lineare-A6}
\left\{\begin{array}{l}
       -\ep (\phi_p (\ep u'))'-C_0 a(x)\phi_p(u)=\lambda \phi_p(u)\\
			\\
			u(s)=0=u(t).
			
			\end{array}
\right.
\end{equation}
Then we have
\begin{equation} \label{num-aut-neg}
\lim_{\ep \to 0^+} \ep l_\ep (s,t)=\dfrac{(C_0)^{1/p}}{\pi_p}\int_s^t \sqrt[p]{a(x)}\,dx.
\end{equation}
\end{lemma}
\vs{8}
\ni
\begin{proof} Let us first recall (see \cite{Zh-01}) that $l_\ep (s,t)$ coincides with the number of zeros in $(s,t)$ of the solution of
\begin{equation} \label{pr11}
\left\{\begin{array}{l}
       \ep (\phi_p (\ep u'))'+C_0 a(x)\phi_p(u)=0\\
			\\
			u(s)=0, \,u'(s)=1.
			\end{array}
\right.
\end{equation}
The equation in \eqref{pr11} is of the form \eqref{eq-riscritta}, with $g(u)=u$ and $C_0a$ instead of $a$; as a consequence, arguing as in the proof of Proposition \ref{stima-zeri}, from \eqref{b14} we immediately deduce the result.
\end{proof}
\vs{12}
\ni
We are now ready to prove Theorem \ref{main-ex-1}, following the same lines of the proof of \cite[Proposition 4.1]{FeMaTa-05}.
\vs{8}
\ni
\begin{proof}[Sketch of proof of Theorem \ref{main-ex-1}] We denote by $\tau$ the vector $(\tau_1,\ldots, \tau_{n^0_\ep+n^1_\ep})$. By compactness, 
\[
\max_{\tau \in \Delta'} f_\ep(\tau)
\]
is attained at a value $\tau^*$. 
\vs{4}
\ni
As a first step, we show that
\begin{equation}\label{nontrivial}
s_0 - h_0 \leq \tau^*_1 < \ldots < \tau^*_{n^0_\ep} \leq t_0 + h_0
\quad \mbox{ and } \quad 
s_1 - h_0 \leq \tau^*_{n^0_\ep + 1} < \ldots < \tau^*_{n^0_\ep + n^1_\ep} \leq t_1 + h_0
\end{equation}
and that the corresponding minimizers $u_{(-)^j} (\ep;\tau^*_j,\tau^*_{j+1})$ are non-zero for every $j=0,\ldots, n^0_\ep+n^1_\ep$.
To this aim, let $\ep>0$ be such that
\begin{equation} \label{nuova1}
n^i_\ep+2<l_\ep (s_i-h_0,t_i+h_0),\quad i=0, 1.
\end{equation}
This choice of $\ep$ is possible since \eqref{ipotesinumerozeri} and \eqref{num-aut-neg} hold and the function $T$ satisfies Proposition \ref{prop-periodo}. Moreover, let $\lambda^i_\ep<0$ be the $(n^i_\ep+2)$-th eigenvalue of \eqref{lineare-A6} in $[s_i-h_0,t_i+h_0]$ and let $e^i_\ep$ be the corresponding eigenfunction, whose zeros we denote by
\[
s_i-h_0=\eta_0^i<\eta_1^i<\ldots<\eta^i_{n^i_\ep+1}<\eta^i_{n^i_\ep+2}=t_i+h_0.
\]
It is trivial to see that there exist $j\in \{0, 1,\ldots, n^0_\ep+1\}$ and $k\in \{0, 1,\ldots, n^0_\ep\}$ such that
\[
[\eta^0_j,\eta^0_{j+1}]\subset [\tau^*_k,\tau^*_{k+1}].
\]
Since the first eigenvalue of \eqref{lineare-A6} in $[\eta^0_j,\eta^0_{j+1}]$ is $\lambda^0_\ep<0$, we deduce that the first eigenvalue of \eqref{lineare-A6} in $[\tau^*_k,\tau^*_{k+1}]$ is
also negative; hence, from Lemma \ref{non-banali} we obtain that
\[
u_{(-)^k} (\ep;\tau^*_k,\tau^*_{k+1})\not\equiv 0.
\]
Hence, $u_\ep \not\equiv 0$. From this, one can show that \eqref{nontrivial} holds true and that \emph{all} the minimizers
$u_{(-)^j} (\ep;\tau^*_j,\tau^*_{j+1})$ are non-trivial just by using the formulas for
the derivatives of $m_{(-)^j}$ contained in Lemmas \ref{minimidir} and \ref{minimidir2}.
We omit the details which can be found in \cite{FeMaTa-05}.
\vs{4}
\ni
\vs{4}
\ni
As a second step, we show that  
\[
s_0-h_0<\tau^*_1,\quad \tau^*_{n^0_\ep}<t_0+h_0,\quad s_1-h_0<\tau^*_{n^0_\ep+1},\quad \tau^*_{n^0_\ep+n^1_\ep}<t_1+h_0.
\]
For instance we check the validity of the relation
\[
s_0-h_0<\tau^*_1,
\]
for $\ep$ sufficiently small. Arguing by contradiction, assume that there exists $\ep_n\to 0^+$ such that
\[
s_0-h_0=\tau^{*,\ep_n}_1.
\]
Hence, $u_{\ep_n}$ is a solution in $I_0=[s_0-h_0,t_0+h_0]$; the corresponding energy $E_{\ep_n}$ satisfies
\[
E_{\ep_n}\to F,
\]
uniformly in $I_0$, for some function $F$ which satisfies \eqref{eq-profilo} and
\[
\int_{I_0} \dfrac{2\sqrt[p]{a(x)}}{T(F(x))}\,dx=\int_{s_0}^{t_0} \dfrac{2\sqrt[p]{a(x)}}{T(E(x))}\,dx.
\]
Using \eqref{num2} and the properties of $T$ we can conclude, as in \cite{FeMaTa-05}, that
\[
F\equiv E, \quad \mbox{in}\ I_0;
\]
in particular, since $E \equiv 0$ on $[s_0-h,s_0]$,
a slight variant of Proposition \ref{prop2.6} implies that
\[
s_0-h_0=\tau^{*,\ep_n}_1<s_0-\dfrac{1}{2}h_0<\tau^{*,\ep_n}_2
\]
for $n$ large. At this point, Lemma \ref{stimeminimi} (iii) can be applied yielding
\[
\dfrac{\partial f_\ep}{\partial \tau_1}(\tau^*)>0
\]
and thus contradicting the fact that $\tau^*$ is a maximizer.
\vs{4}
\ni
Hence, we have shown that $\tau^*\in \mathring \Delta'$; using again the formulas for
the derivatives of $m_{(-)^j}$ contained in Lemmas \ref{minimidir} and \ref{minimidir2}
(compare, in particular, with \cite[Proposition A8]{FeMaTa-05}) 
this is sufficient to prove that the function $u_\ep$ defined in \eqref{num} is a solution of \eqref{bvp-1} with the required properties.
\end{proof}
\vs{12}
\ni
It is clear that Theorem \ref{main-ex-1} can be extended to the case when the support of $E$ is the union
of finitely many intervals $(s_i,t_i)$ satisfying the non-degeneracy condition \eqref{num2}.
Using an approximation argument developed in \cite{FeMaTa-05}, the general case can be treated as well.
Summing up, we can finally state the following existence result:
\vs{12}
\ni
\begin{theorem} \label{main-ex-2} For every solution $E$ of \eqref{eq-profilo} there exists a family of solutions $u_\ep$ of \eqref{bvp-1} such that
\[
E_\ep \to E,\quad \mbox{uniformly in $[0,1]$}.
\]
\end{theorem}
\vs{12}
\ni
{\bf Acknowledgement.} The authors thank Prof. Susanna Terracini for having suggested the topic and for useful discussions.

\vspace{1 cm}
\normalsize

Authors' addresses:
\bigbreak
\medbreak
\indent Alberto Boscaggin \\
\indent Dipartimento di Matematica e Applicazioni, Universit\`a di Milano Bicocca, \\
\indent Via Cozzi 53, I-20125 Milano, Italy \\
\indent e-mail: alberto.boscaggin@unimib.it \\

\medbreak
\indent Walter Dambrosio \\
\indent Dipartimento di Matematica, Universit\`a di Torino \\
\indent Via Carlo Alberto 10, I-10123 Torino, Italy \\
\indent e-mail: walter.dambrosio@unito.it

\end{document}